\documentclass[11pt,reqno]{amsart}
\usepackage{amsmath,amssymb,amsthm,amsfonts,enumerate,hyperref}
\usepackage{graphicx}
\theoremstyle{plain}
\newtheorem{theorem}{Theorem}[section]
\newtheorem{proposition}[theorem]{Proposition}
\newtheorem{cor}[theorem]{Corollary}
\newtheorem{lemma}[theorem]{Lemma}

\theoremstyle{definition}
\newtheorem{definition}[theorem]{Definition}

\theoremstyle{remark}

\numberwithin{equation}{section}

\begin{document}

\title{Computing Perelman's $\nu$-functional}

\author{Stuart James Hall}

\address{Department of Mathematics, Imperial College, London SW7 2AZ, UK\\
email: stuart.hall06@imperial.ac.uk}

\begin{abstract}
This is a short note in which we show how to calculate the value of Perelman's $\nu$-functional for a variety of metrics. In particular we complete the calculation of values for the known 4-dimensional Einstein and shrinking Ricci soliton metrics.
\end{abstract}
\maketitle

\section{Introduction}
\subsection{Overview}
One of Perelman's breakthroughs in his series of preprints concerning the Ricci flow was a functional that behaved a lot like an energy functional for the flow, the $\nu$-functional \cite{Per}. This functional has some rather desirable qualities; its critical points are Ricci solitons (a class of Riemannian metrics generalising Einstein metrics) and it is increasing along a Ricci flow and stationary only at Ricci solitons. Hence, on manifolds admitting more than one soliton, it is extremely interesting to know the precise value of the function at critical points in order to determine how the geometry of the manifold evolves under the Ricci flow. This article computes the values of the $\nu$-functional on some $4$-dimensional manifolds that are known to admit two Ricci solitons (one being an Einstein metric, the other not). We also show how to compute the value of the $\nu$-functional on a class of solitons called toric-K\"ahler-Ricci solitons, the existence of which is due to Wang and Zhu \cite{WangZhu}. In particular we complete the calculation for all known solitons on compact 4-dimensional manifolds, a task begun by Cao-Hamilton-Ilmanen \cite{CHI}.\\

The main focus of our results will be on the Fano surfaces $\mathbb{CP}^{2}\sharp\overline{\mathbb{CP}}^{2}$ and $\mathbb{CP}^{2}\sharp2\overline{\mathbb{CP}}^{2}$. Each of these manifolds admits a Hermitian, non-K\"ahler, Einstein metric. The existence of these metrics is due to Page \cite{Page} in the $\mathbb{CP}^{2}\sharp\overline{\mathbb{CP}}^{2}$ case and Chen-LeBrun-Weber \cite{CLB} in the $\mathbb{CP}^{2}\sharp2\overline{\mathbb{CP}}^{2}$ case. Due to the work of Derdzinski \cite{Derd} it is known that both of these metrics are conformally K\"ahler.   We will find it convenient to compute the quantity $\Theta(g)=e^{\nu(g)}$; $\Theta(g)$ is known as the Gaussian density of $g$ and was defined in Cao-Hamilton-Ilmanen \cite{CHI}. There is a general formula, given by Cao-Hamilton-Ilmanen in \cite{CHI},  for the Gaussian density of an Einstein metric with positive scalar curvature which gives the density in terms of the volume and the scalar curvature of the metric (\emph{cf}. Proposition \ref{P1}). For the Page and Chen-LeBrun-Weber metric we can compute the density in terms of topological data. The value of the $\nu$-functional in each case, follows from Theorem \ref{T1}. 
\begin{theorem}\label{T1}
Let $(M^{4},g,J)$ be a K\"ahler manifold such that the K\"ahler metric $g$ is conformal to an Einstein metric $g_{E}$. 
The Gaussian density of $g_{E}$ is given by
$$\Theta(g_{E})=\frac{3}{2e^{2}}(2\chi(M)+3\sigma(M))-\frac{2Cal(g)}{(8\pi e)^{2}}.$$
Here $\chi(M)$ is the Euler characteristic of $M$, $\sigma(M)$ is the signature of $M$ and $Cal(g)$ is the Calabi energy of the K\"ahler metric $g$,
$$Cal(g)=\int_{M}R(g)^{2}dV_{g}$$
where $R(g)$ is the scalar curvature of $g$.
\end{theorem}

Both of the Fano surfaces $\mathbb{CP}^{2}\sharp\overline{\mathbb{CP}}^{2}$ and $\mathbb{CP}^{2}\sharp2\overline{\mathbb{CP}}^{2}$ belong to a class of complex manifolds called toric-K\"ahler manifolds. In fact, the manifold $\mathbb{CP}^{2}\sharp\overline{\mathbb{CP}}^{2}$ belongs to an even more special class of manifolds as it admits a cohomogeneity one action by the group $U(2)$. Wang and Zhu showed the existence of a unique (up to automorphisms) K\"ahler-Ricci soliton on any toric-K\"ahler manifold and, on $\mathbb{CP}^{2}\sharp\overline{\mathbb{CP}}^{2}$ and $\mathbb{CP}^{2}\sharp2\overline{\mathbb{CP}}^{2}$, the solitons are not Einstein (non-trivial). The soliton on  $\mathbb{CP}^{2}\sharp\overline{\mathbb{CP}}^{2}$ is $U(2)$-invariant and had previously been constructed independently by Koiso \cite{Koiso} and Cao \cite{Cao}.  Toric-K\"ahler manifolds have a rich theory that we shall explain a little more in section 3. For now we simply note our result that the value of the Gaussian density for these solitons follows from:
\begin{theorem}\label{T2}
Let $M^{2n}$ be a toric-Fano manifold and let $g_{KRS}$ be the  Wang-Zhu K\"ahler-Ricci soliton normalised so that
$Ric(g_{KRS})+Hess(f)=g_{KRS}$ for a smooth potential function $f$. Let $P\subset \mathbb{R}^{n}$ be the moment polytope translated so that center of mass is at $0$. If we define the convex function $F:\mathbb{R}^{n}\rightarrow \mathbb{R}$ by
$$F(\underbar{c})=\int_{P}e^{-c\cdot x}dx,$$
the Gaussian density of the soliton is given by
$$\Theta(g_{KRS}) = e^{-n}\min_{\mathbb{R}^{n}}F.$$ 
\end{theorem}
As the Page metric and Koiso-Cao soliton on $\mathbb{CP}^{2}\sharp\overline{\mathbb{CP}}^{2}$ are explicit, the value of the $\nu$-functional at these metrics was already calculated \cite{CHI}.  The Chen-LeBrun-Weber metric and the Wang-Zhu soliton on $\mathbb{CP}^{2}\sharp2\overline{\mathbb{CP}}^{2}$ are not constructed explicitly and so computing the exact value of the $\nu$-functional is more challenging.  Recently Headrick and Wiseman \cite{HW} used a numerical approximation to the Wang-Zhu soliton to compute its $\nu$-value.  Our methods are explicit and very simple once one has the huge machinary of toric-K\"ahler manifolds at one's disposal. The main contribution of this paper is the value of $\nu$ at the Chen-LeBrun-Weber metric. It seems a little unexpected that this value is greater than that of the Wang-Zhu soliton \cite{Caocomm}.

\subsection{Ricci solitons and the $\nu$-functional}
We begin by recalling the notion of Ricci soliton.
\begin{definition}[Ricci soliton] A Riemannian metric g is called
a Ricci soliton if there exists a vector field V and constant \(\rho\) satisfying \begin{displaymath}
Ric(g)+L_{V}g=\rho g.  
\end{displaymath}
Moreover, if V is the gradient of a smooth function \(f\), the  soliton is
called a gradient Ricci soliton with potential function \(f\). For \(\rho > 0\) the soliton is called shrinking,
for \(\rho = 0\) the soliton is called steady and for \(\rho < 0\) the soliton is called expanding.    
\end{definition}
We see that by setting \(V=0\) in the above definition we recover the notion
of an Einstein metric so solitons can be thought of as generalisations of Einstein metrics. We shall henceforth refer to Einstein metrics as being trivial Ricci solitons. In this article we will only be concerned with shrinking Ricci solitons on compact manifolds. The following theorem of Perelman means we can restrict to gradient Ricci solitons.
\begin{theorem}[Perelman \cite{Per}]
Let $g$ be a shrinking Ricci soliton on a compact manifold $M$, then $g$ is a gradient Ricci soliton.
\end{theorem}

 If the metric $g$ is K\"ahler then the vector field  $V$ is holomorphic and we say that $g$ is a K\"ahler-Ricci soliton.  It is straightforward to see that shrinking K\"ahler-Ricci solitons can only occur on Fano manifolds.

Perelman defined the following functionals:
\begin{definition}[\(\mathcal{W}\)-functional] If \((M^{n},g)\) is a closed Riemannian
manifold, \(f\) a smooth function on \(M\) and \(\tau >0\) a real number
then we define
\begin{displaymath}
\mathcal{W}(g,f,\tau) = \int_{M}[\tau(R+|\nabla f|^{2})+f-n](4\pi\tau)^{-\frac{n}{2}}e^{-f}dV_{g}
\end{displaymath}
where \(R\) is the scalar curvature of \(g\).  
\end{definition}
\begin{definition}[\(\nu\)-energy \& Gaussian density]
Let \(M^{n}\), \(g\), \(f\) and \(\tau\) be as above, then
\begin{displaymath}
\nu(g)=\inf \{\mathcal{W}(g,f,\tau) : \frac{1}{(4\pi\tau)^{\frac{n}{2}}}\int_{M}e^{-f}dV=1\}.
\end{displaymath}
 The Gaussian density of \(g\),
 \(\Theta(g)\), is \begin{displaymath}
\Theta(g) = e^{\nu(g)}.
\end{displaymath}
\end{definition}
The only comment on existence theory we will make is that for a fixed $g$ there exists a smooth $f$ and $\tau>0$ that achieve the infimum in the definition of the $\nu$-functional. The minimising pair $(f,\tau)$ solve the equations
$$\tau(-2\Delta f +|\nabla f|^{2}-R)-f+n+\nu =0 $$
and
$$\frac{1}{(4\pi\tau)^{n/2}}\int_{M}fe^{-f}dV_{g}=\frac{n}{2}+\nu(g). $$ 

The main theorem concerning the $\nu$-functional is the following:
\begin{theorem}[Perelman \cite{Per}]
Let $(M,g)$ be a closed Riemannian manifold. If $g$ is a critical point of the $\nu$-functional then $g$ is a gradient Ricci soliton and the compatible triple $(g,f,\tau)$ satisfy
$$Ric(g)+Hess(f)=\frac{1}{2\tau}g.$$
\end{theorem}
\section{Calculation for conformally K\"ahler, Einstein metrics}
The purpose of this section is to prove Theorem 1. If the soliton is trivial then there is the following well-known formula for the Gaussian density.
\begin{proposition}\label{P1}[\cite{CHI} for example]
 If \((M^{n},g)\) is an Einstein
manifold of positive scalar curvature then
\begin{displaymath}
\Theta(g)=\left(\frac{R(g)}{2\pi ne}\right)^{\frac{n}{2}}Vol(M,g)
\end{displaymath}
where \(R(g)\) is the scalar curvature of \(g\) and \(Vol(M,g)\) is the volume
of \(M\) with respect to \(g\).
\end{proposition}
As mentioned in the Introduction, the Einstein metrics we are interested in are conformally K\"ahler.  We will need the following well-known fact about the curvature tensor of K\"ahler metrics:
\begin{lemma}\label{L1}[\cite{Derd} for example]
Let $(M,g,J)$ be a K\"ahler manifold. Then
$$|W^{+}(g)|^{2} = \frac{R^{2}(g)}{24}$$
where $W^{+}(g)$ is the self-dual part of the Weyl tensor of $g$ and $R(g)$ is the scalar curvature of $g$.
\end{lemma}
The proof of Theorem 1 is now in our hands
\begin{proof}(Of Theorem \ref{T1}.)
We let $g_{E}$ be an Einstein metric conformal to a K\"ahler metric $g$. We have the Allendorfer-Weil formula 
$$\chi(M)=\frac{1}{8\pi^{2}}\int_{M}|W(g_{E})|^{2}+\frac{R^{2}(g_{E})}{24}dV_{g_{E}} $$
where $W(g_{E})$ and $R(g_{E})$ are the Weyl tensor and  scalar curvature of $g_{E}$ respectively. As the scalar curvature is constant we have
$$R(g_{E})^{2}Vol(M,g_{E})=192\pi^{2}\chi(M)-24\int_{M}|W(g_{E})|^{2}dV_{g_{E}}.$$
The integral
$$\mathcal{F}=\int_{M}|W(g_{E})|^{2}dV_{g_{E}}$$
is a conformal invariant so we can compute it with respect to the K\"ahler metric $g$. The Hirzebruch signature  formula,
\begin{displaymath}
\sigma(M)=\frac{1}{12\pi^{2}}\int_{M}|W^{+}(g)|^{2}-|W^{-}(g)|^{2}dV_{g},
\end{displaymath}  
allows use to write \(\mathcal{F}\) as
\begin{displaymath}
\mathcal{F}= \int_{M}2|W^{+}(g)|^{2}dV_{g}-12\pi^{2}\sigma(M)
\end{displaymath}
where \(W^{+}(g)\) and \(W^{-}(g)\) denote the self-dual and anti-self-dual parts
of the Weyl tensor of $g$. Hence by Lemma \ref{L1} we have
$$\mathcal{F} = \frac{Cal(g)}{12}-12\pi^{2}\sigma(M).$$
Substituting these formulae into Proposition \ref{P1} the result follows.
\end{proof}
The calculation of $\nu$ for the Page and Chen-LeBrun-Weber metrics reduces to computing the Calabi energy of the K\"ahler metrics that they are conformally equivalent to. The K\"ahler metrics are themselves rather special; they are extremal metrics meaning that they represent critical points of the Calabi energy when one varies through K\"ahler metrics in the same K\"ahler class. Due to their distinguished status, one can compute the Calabi energy of these K\"ahler metrics without knowing the metric explicitly. We shall simply state the result of this calculation.
\begin{cor}
The Gaussian density of the Page metric is 0.5172 to 4 decimal places. The Gaussian density of the Chen-LeBrun-Weber metric is 0.4552 to 4 decimal places.
\end{cor}
\begin{proof}
The Euler characteristic of $\mathbb{CP}^{2}\sharp\overline{\mathbb{CP}}^{2}$ is 4 and its signature is 0. In \cite{LeBrun}, LeBrun computes that the Calabi energy of the extremal K\"ahler metric conformal to the Page metric is given 
by \(96\pi^{2}\min_{\mathbb{R}}h\), where \(h\) is the function
\begin{displaymath}
h(x)=\left(\frac{4+14x+16x^{2}+3x^{3}}{x(6+6x+x^{2})}\right).
\end{displaymath}
The function \(h\) is minimised when \(x \approx 2.183933\) giving a value
of 2.72621; putting these
values into Theorem \ref{T1} we obtain that the density of the Page metric is
0.5172 to 4 decimal places.\\
\\
 The Euler characteristic of $\mathbb{CP}^{2}\sharp2\overline{\mathbb{CP}}^{2}$ is 5 and its signature is -1. In \cite{CLB}, Chen, LeBrun and Weber compute that the Calabi energy of the extremal K\"ahler metric conformal to the Chen-LeBrun-Weber metric is given 
by \(32\pi^{2}\min_{\mathbb{R}}f\), where \(f\) is the function
\begin{displaymath}
f(x)=3\left(\frac{32 + 176x +318x^2 +280x^{3}+132x^{4}+32x^{5}+3x^{6}}{12+72x+
138x^{2}+120x^{3}+54x^{4}+12x^{5} +x^{6}}\right).
\end{displaymath}
The function \(f\) is minimised when \(x \approx 0.95771\) giving a value
of 7.13647; putting these
values into Theorem \ref{T1} we obtain that the density of the Chen-LeBrun-Weber metric is
0.4552 to 4 decimal places.
\end{proof}
\section{Calculation for K\"ahler-Ricci solitons}
A gradient Ricci soliton is given by a triple $(g,f,\tau)$ such that
$$Ric(g)+Hess(f)=\frac{1}{2\tau}g.$$ 
The triple satisfy the compatibility condition
$$\frac{1}{(4\pi\tau)^{n/2}}\int_{M}e^{-f}dV_{g}=1$$ and, by examining the Euler-Lagrange system for the triple, they also satisfy the equation
$$\frac{1}{(4\pi\tau)^{n/2}}\int_{M}fe^{-f}dV_{g}=\frac{n}{2}+\nu(g).$$
If we scale the soliton so that $\tau=0.5$ and then define the quantity
$$Z(\beta) = \frac{1}{(2\pi e)^{n/2}}\int_{M}e^{-\beta f}dV_{g},$$
the value of $\nu(g)$ is given by
$$\nu(g)=(1-\beta\frac{d}{d\beta})\log Z(\beta)\bigg|_{\beta =1}.$$
\subsection{Toric-K\"ahler-Ricci solitons}
In this section we discuss some of the theory of toric-K\"ahler metrics as applied to Ricci solitons. In the interests of brevity we will not go into detail but an interested reader should consult the recent survey by Donaldson \cite{DonTor}. The main feature of the theory we will use is that, for a Fano toric-K\"ahler manifold $M^{2n}$, there is a dense open subset $M^{\circ}\subset M$ such that $M^{\circ} = P \times \mathbb{T}^{n}$ where $P\subset\mathbb{R}^{n}$ is a reflexive polytope with center of mass at the origin. The polytope $P$ is called the \emph{moment polytope}. In these coordinates the metric $g$ is given by
$$g = u_{ij}dx_{i}dx_{j}+u^{ij}d\theta_{i}d\theta_{j} $$ 
where $u_{ij}$ is the Euclidean Hessian of a special convex function $u:P\rightarrow \mathbb{R}$ and $u^{ij}$ is the matrix inverse of $u_{ij}$. This means that the volume form in these coordinates is given by $dx_{1}\wedge dx_{2}\wedge ...\wedge dx_{n} \wedge d\theta_{1} \wedge d\theta_{2} \wedge ...\wedge d\theta_{n}$.
\begin{lemma}\label{TorKRS}[\cite{DonTor} for example]
Let $(g,f,\tau)$ be a toric-K\"ahler-Ricci soliton with reflexive moment polytope $P$ with a center of mass $0$. If we scale so that the soliton satisfies
$$Ric(g)+Hess(f)=g$$ then in these coordinates the soliton potential function is given by $f(x) = a_{1}x_{1}+...+a_{n}x_{n}$. Furthermore the coefficients $a_{i}$ are determined by the constraint
$$\int_{P}x_{i}e^{-f}dx_{1} \wedge .. \wedge dx_{n} = 0.$$
\end{lemma}
We can now prove Theorem \ref{T2}
\begin{proof}(Of Theorem \ref{T2}.)
The constraint on $f$ in Lemma \ref{TorKRS} can be viewed as saying that the coefficients $(a_{1},...,a_{n})$ are critical points of the function $F:\mathbb{R}^{n}\rightarrow \mathbb{R}$ given by
$$ F(c)= \int_{P}e^{-c\cdot x}dx.$$ It is straightforward to see that this function is convex and proper and so $(a_{1},...,a_{n})$ is the point where it achieves its minima. The $\nu$-energy is given by 
$$\nu(g_{KRS})=(1-\beta\frac{d}{d\beta})\log Z(\beta)\bigg|_{\beta =1}.$$  
The function $Z(\beta)$ is
$$\frac{1}{(2\pi e)^{n}}\int_{M}e^{-\beta f}dV_{g} = \frac{1}{(2\pi e)^{n}}\int_{P\times \mathbb{T}^{n}}e^{-\beta f}dx\wedge d \theta = \frac{1}{e^{n}}\int_{P}e^{-\beta f}dx,$$ 
where the last equality follows from the fact that the volume of $\mathbb{T}^{n}=(2\pi)^{n}$. The constrain in Lemma \ref{TorKRS} means that $\frac{d}{d\beta}\log Z(\beta) |_{\beta =1}=0$ and so
$$\nu(g_{KRS}) = \log Z(1).$$
Hence $\Theta(g_{KRS}) = e^{-n}\min_{\mathbb{R}^{n}}F(c)$.
\end{proof}
\begin{cor}
The Gaussian density of the Koiso-Cao soliton is 0.5179 to 4 decimal places.  The Gaussian density of the Wang-Zhu soliton is 0.4549 to 4 decimal places.
\end{cor}
\begin{proof}
For the Koiso-Cao soliton on \(\mathbb{CP}^{2}\# \overline{\mathbb{CP}^{2}}\), the moment polytope is the
trapezium \(T\) with vertices at (2,-1), (-1,2), (-1,0)
and (0,-1). We can take the potential function to be invaraint under the bilateral symmetry $x_{1} \leftrightarrow x_{2}$ so  the soliton potential function \(f\) is of the form
\begin{displaymath}
f=c(x_1+x_2)
\end{displaymath} 
where the constant \(c\) is the unique minimiser of the function  
 \begin{displaymath}
F(c)= \int_{T}e^{-c(x_{1}+x_{2})}dx.
\end{displaymath}
Evaluating, we have that
\begin{displaymath}
F(c)=\frac{(e^{2c}-e^{-c}-3ce^{-c})}{c^{2}}.
\end{displaymath}
This function is minimised at \(c \approx 0.5276\) giving a value of 3.8266.
 The density of the Koiso-Cao soliton \(g_{KC}\) is
 \begin{displaymath}
\Theta(g_{KC})=\frac{\min_{\mathbb{R}}H}{e^{2}}=\frac{3.8266}{e^{2}}=0.5179
\end{displaymath}
to 4 decimal places.

 For the Wang-Zhu soliton on \(\mathbb{CP}^{2}
\#2 \overline{\mathbb{CP}^{2}}\) the moment polytope is the pentagon $P$ with vertices at (-1,-1),
(1,-1), (1,0), (0,1) and (-1,1). Again the potential function is invariant under the symmetry $x_{1} \leftrightarrow x_{2}$  so the soliton potential function \(f\) is of the form
\begin{displaymath}
f=c(x_1+x_2)
\end{displaymath} 
where the constant \(c\) is the unique minimiser of the function  
 \begin{displaymath}
F(c)= \int_{P}e^{-c(x_{1}+x_{2})}dx.
\end{displaymath}
Evaluating, we have that
\begin{displaymath}
F(c)=\frac{(e^{2c}-2+(1-c)e^{-c})}{c^{2}}.
\end{displaymath}
This function is minimised at \(c \approx -0.434748\) giving a value of 3.36094.
 The density of the Wang-Zhu soliton \(g_{WZ}\) is
 \begin{displaymath}
\Theta(g_{WZ})=\frac{\min_{\mathbb{R}}H}{e^{2}}=\frac{3.36094}{e^{2}}=0.4549
\end{displaymath}
to 4 decimal places.
\end{proof} 
\emph{Acknowledgements:}  This work forms part of the author's Ph.D thesis at Imperial College, London funded by the Engineering and Physical Sciences Research Council.  I would like to thank my advisor Simon Donaldson for his useful comments and encouragement during the course of this work.  I would also like to thank Professor Huai-Dong Cao and Thomas Murphy for comments on earlier versions of this article. I would also like to thank the anonymous referee for helpful comments.

\end{document}